\documentclass[a4paper,12pt]{article}

\usepackage{pstricks}
\usepackage{graphicx,psfrag}
\usepackage{amsmath}
\usepackage{amssymb}
\usepackage{amsthm}
\usepackage{color}

      \newcommand {\al}   {\alpha}          \newcommand {\bt}  {\beta}
      \newcommand {\gam } {\gamma}

                 \newcommand {\vphi} {\varphi}
      \newcommand {\lam}  {\lambda}

      \newcommand {\pl}   {\partial}        \newcommand {\s}    {\sigma}

             \newcommand {\RRR}  {{\mathbb R}}

     \newcommand {\beq}  {\begin{equation}}
      \newcommand {\eeq}  {\end{equation}}

      \newtheorem{theorem}{Theorem}
      \newtheorem{lem}{Lemma}
      
      \newtheorem{zam}{Remark}
      \newtheorem{opr}{Definition}

\author{Alexander Plakhov\thanks{University of Aveiro, Portugal}{ }and Vera Roshchina\thanks{Collaborative Research Network, University of Ballarat, Australia}}

\title{Bodies with mirror surface invisible from two points}

\date{}

\begin{document}

\maketitle

\begin{abstract}
Here we are concerned with a special issue of billiard invisibility, where a bounded set with a piecewise smooth boundary in Euclidean space is identified with a body with mirror surface, and the billiard in the complement of the set is identified with the dynamics of light rays outside the body in the framework of geometric optics. We show that in this setting it is possible to construct a body invisible from two points.
\end{abstract}



\section{Introduction}
The problems studied in the framework of billiard invisibility involve the mathematical design of bodies with well-defined surfaces whose scattering map preserves certain trajectories of a flow of elastic particles. The main practical application of this study is the optical shielding: by surrounding an object by a specially designed mirror surface, it is possible to create an illusion of invisibility from given points or directions.

From the practical viewpoint, this approach is a low-tech alternative to the major modern attempts at achieving invisibility, which are primarily focussed around the design of metamaterials that allow the bending of electromagnetic waves around the concealed object. The idea was first suggested in \cite{PendryEtAl06}, and tangible results in this direction were achieved by research teams at Duke University, first in 2006 when a imperfect prototype of microwave-range cloaking device was created \cite{SMJCPSS2006}, and then in 2012, when remarkable results were achieved completely cloaking a centimetre-thin board, even though in one direction only. Another development in this direction \cite{ZhangEtAl2011} is a successful concealing of a nanoscale object under a carpet cloak made of layers of silicon oxide and silicon nitride arranged in a special way.

Another approach to invisibility was recently suggested in \cite{ZhangBaiLeEtAl2011} and lead to the development of a working invisibility cloak based on calcite crystals (also see \cite{ChenEtAl2011}). The cloak works only under one light polarization: it is essentially two-dimensional, although works at all angles.

The first work that targets the problem of designing a body invisible in a direction in the framework of mirror invisibility appears in \cite{AlekPla2009} and is motivated by the problem of constructing a nonconvex body or zero resistance. The authors demonstrated that there exists a (connected and even simply connected) body invisible in one direction: if this body is manufactured out of perfectly reflective mirrors, a laser beam sent through this construction in the direction of invisibility would leave the body along the same trajectory. Remarkably, in \cite{Lak2012} it is shown that this body is also invisible in acoustic waves.

This pioneering research led to several intriguing mathematical problems. Some of them, proposed by Sergei Tabachnikov \cite{Tabach}, ask whether it is possible to design a body with a mirror surface invisible in two directions or a body invisible from a point. The former problem was solved in \cite{PR2011}: it was demonstrated that a parabolic construction can be used to produce a body invisible in two directions in three-dimensional case. This body consists of two connected components. It is possible to construct two- and three-dimensional infinitely connected fractal bodies invisible in two and three directions respectively \cite{PR2013}. Note that it is impossible to construct a body invisible in all directions \cite{PR2011}, and at least in two-dimensional case there is no piecewise smooth body invisible in a countable number of directions \cite{PCInvisImposs}. A somewhat related development that should be mentioned here is the phenomena best known under the term of digital sundial. It is possible to construct fractal bodies such that their projections on almost all planes can be prescribed up to a set of zero measure (see \cite{Falc2003}). In this framework, Burdzy and Kulczycki \cite{BurKul2012} showed that there are two-dimensional bodies which consist of an infinite number of linear segments, and which are almost invisible in almost all directions.

In addition to the fascinating research problems discussed above, billiard invisibility is linked to some fundamental properties of dynamical systems. As it is pointed out in \cite{ebook} (see also \cite{Gut2012}), billiard invisibility is related to the long-standing Ivrii's conjecture. The Ivrii's conjecture states that the measure of periodic trajectories in a piecewise smooth planar billiard is zero, and, in particular, for any natural number $n$ the measure of periodic trajectories of period n is zero. The conjecture was proved for several important special cases (the most recent breakthrough is \cite{GluKudr2012}, where the authors proved the nonexistence of billiards with an open set of periodic quadrilateral orbits).

Recently we have identified bodies (connected but not simply connected) invisible from a point \cite{ebook, invis1point}. In this paper we provide a construction of a body invisible from 2 points. The body has infinitely many connected components. To our best knowledge, this construction demonstrates a new property of confocal conics. Thus, if you look with two eyes open from a fixed position, the body disappears --- becomes completely invisible. The previous result meant that you needed to look at the body with one eye open and the other closed.

The main result of the paper is the following Theorem \ref{t1}.

\begin{theorem}\label{t1}
Given two different points in $\RRR^d$, $n = 2,\, 3,\ldots$, there exists a body in $\RRR^d$ invisible from these points.
\end{theorem}

We also show that mixed-type invisibility, where two points are substituted with a point and a direction, also holds true.

\begin{theorem}\label{t2}
Given a points $A \in \RRR^d$ and a vector $v \in S^{d-1}$, there exists a body in $\RRR^d$ invisible simultaneously from the point $A$ and in the direction $v$.
\end{theorem}

The rest of the paper is devoted to the proof of these results.

\section{Invisibility from two points}

We begin with reminding the relevant definitions, then explain our construction and prove that it is invisible from two points.

\begin{opr}\label{def3}\rm
A {\it body} is a finite or countable union of its connected components, where each component is either a bounded domain with a piecewise smooth boundary, or a bounded piece of smooth hyper-surface.
\end{opr}

\begin{opr}\label{def1}\rm
A body $B \subset \RRR^d$ is said to be {\it invisible from a point} $O \in \RRR^d \setminus B$, if for almost all $v \in S^{d-1}$ the billiard particle in $\RRR^d \setminus B$ emanating from $O$ with the initial velocity $v$, after a finite number of reflections from $\pl B$ will eventually move freely with the same velocity $v$ along a straight line containing $O$.
\end{opr}

If the point $O$ is infinitely distant, we get the notion of a body invisible in a direction. Namely, we have the following definition.

\begin{opr}\label{def2}\rm
A body $B$ is {\it invisible in a direction} $v \in S^{d-1}$, if for almost all straight lines in $\RRR^d$ with the director vector $v$, the billiard particle in $\RRR^d \setminus B$ that initially moves along this line with velocity $v$, after a finite number of reflections from $\pl B$ will finally move freely along the same line with the same velocity $v$.
\end{opr}

\subsection{Construction of a body invisible from two points}

In this section we construct a 2D body invisible from two different points $A_1$ and $A_2$. An invisible body in three or higher dimensions is obtained by rotating the 2D body around the axis $A_1A_2$.

\begin{figure}[h]
\centering
\includegraphics[scale=1]{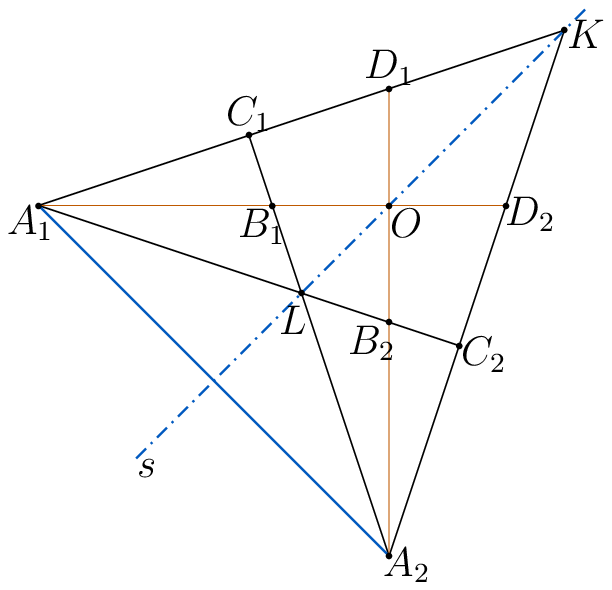}\qquad
\includegraphics[scale=1]{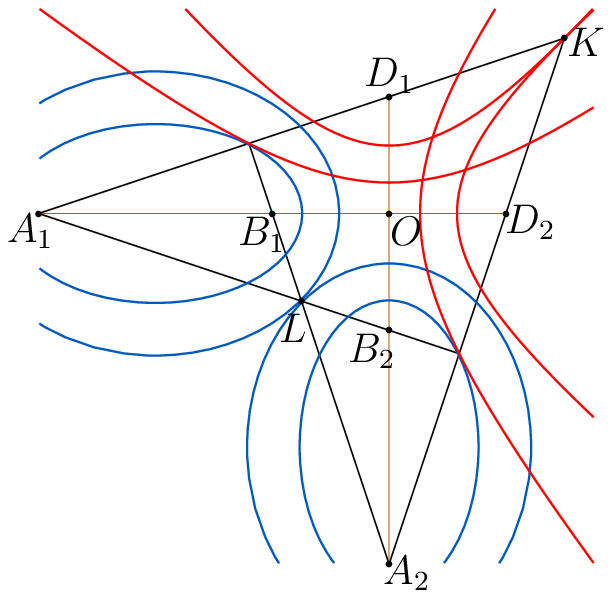}\\
(a)\hskip 180pt (b)
\caption{Construction of the eight basic curves}
\label{fig:v01}
\end{figure}

Let $s$ be the symmetry axis of the system $\{A_1, A_2\}$, i.e. $s$ is a straight line perpendicular to the segment $A_1 A_2$.
Choose two different arbitrary points on $s$ such that they both lie on the same side of the segment $A_1A_2$ (see Fig.~\ref{fig:v01} (a)).
Denote the point nearest to $A_1A_2$ by $L$ and the other one by $K$. Let $C_1$ be the intersection of  $A_1K$ with the line through $A_2$ and $L$, and let $C_2$ be the symmetric intersection of $A_2K$  with the line through $A_1$ and $L$ (see Fig.~\ref{fig:v01} (a)). Choose a point $O$ on $LK$. Let $D_2$ be the intersection of $A_2 K$ with the line through $A_1$and $O$ and let $B_1$ be the intersection of the segments $A_1D_2$ and $C_1L$. The points $D_1$ and $B_2$ are constructed symmetrically (see Fig.~\ref{fig:v01} (b)).

Draw two confocal ellipses with foci $A_1$ and $B_1$ through the points $C_1$ and $L$ and two confocal hyperbolas with foci $A_2$ and $D_1$ through the points $C_1$ and $K$. Repeat this symmetrically on the other side of $s$ (see Fig.~\ref{fig:v01} (b)).
Thus, we have 8 curves of second order through the points $C_1,\, K,\, C_2,\, L$; each curve is either an ellipse or a hyperbola; there are 2 curves through each point. In the sequel by hyperbola we mean the branch of hyperbola that passes through the corresponding point. Each curve (an ellipse or a branch of hyperbola) bounds a convex set; by exterior of an ellipse or a hyperbola we mean the exterior of the corresponding set.

Observe that the curves through each of the points $C_1,\, K,\, C_2,\, L$ are tangent to each other at these points, and hence intersect only at these points. This is due to the property shared by ellipses and parabolas: the tangent line at a point on the curve always makes equal angles with the lines drawn through this point and foci. Hence it is not difficult to see that all of the 8 curves bisect the relevant angles.

Take $H_1$ on the bisector of the angle $B_1C_1D_1$ inside the quadrangle $B_1C_1D_1O$ and denote by $H_2$ the point symmetric to $H_1$ with respect to $KL$. Further, let $N$ be the point of intersection of $A_1H_2$ with $A_2H_1$, and let $M$ be the point of intersection of $A_1H_1$ with $A_2H_2$. Obviously, $N$ belongs to the quadrangle $OB_2LB_1$ and $M$ belongs to the quadrangle $OD_1KD_2$, and both points belong to the line $KL$, which is the bisector of the angles $B_1LB_2$ and $D_1KD_2$. Since the chosen points lie on the bisectors, they belong to the exterior of the corresponding pairs of curves.

If the symmetric configuration of points is slightly changed, one can also choose the points $H_1,\, M,\, H_2,\, N$ possessing the mentioned properties.

\begin{figure}[h]
\centering
\includegraphics[scale=0.8]{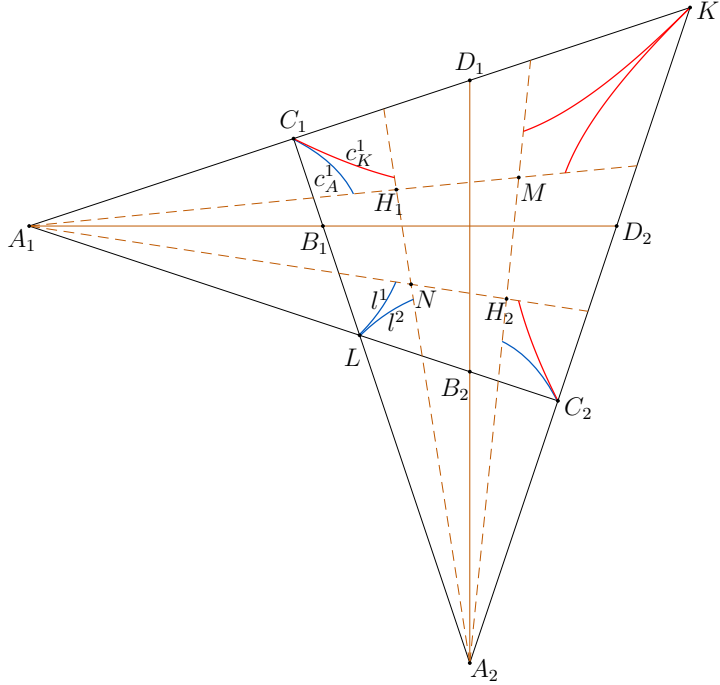}\qquad
\caption{Construction of the initial curvilinear segments}
\label{fig1}
\end{figure}

Now, let $l^1$ be the arc of the ellipse with foci $A_1$ and $B_1$ with the endpoints at $L$ and at the point of intersection of the ellipse with the line $A_1N$, and $l^2$ be the arc of the ellipse with foci $A_2$ and $B_2$ with the endpoints at $L$ and at the point of intersection of the ellipse with the line $A_2N$. It follows from the construction that these points of intersection lie on the segments $A_1N$ and $A_2N$, respectively, and the arcs $l^1$ and $l^2$ belong to the quadrangle bounded by the straight lines $A_1N$,\, $A_2N$,\, $C_1L$,\, $C_2L$. Let $Q_L$ be the curvilinear quadrangle bounded by $l_1$,\, $l_2$, and by segments of the lines $A_1N$,\, $A_2N$. In a similar way define the curvilinear quadrangles $Q_{C_1}$,\, $Q_{C_2}$, and $Q_{K}$; each of them is bounded by two segments of lines $A_1M$,\, $A_2M$,\, $A_1N$,\, $A_2N$ and by two arcs. There are 8 arcs in total (four of them: $l^1$,\, $l^2$,\, $c_A^1$,\, $c_K^1$, are shown in Fig.~\ref{fig1});
each of them generates an infinite sequence of arcs to be defined below, 8 sequences in total. Each sequence is totally contained in one of 8 quadrangles that are shown in figure and determined by their diagonals: $NB_1$,\, $NB_2$,\, $B_1H_1$,\, $B_2H_2$,\, $H_1D_1$,\, $H_2D_2$,\, $MD_1$,\, $MD_2$. Let us describe the sequences $l_0,\, l_1,\, l_2, \ldots$ and $c_0,\, c_1,\, c_2, \ldots$ generated by the arcs $l^2$ and $c_K^1$, respectively; the other sequences are defined analogously.

The sequence $l_0,\, l_1,\, l_2, \ldots$ of arcs of ellipses with foci at $A_2$ and $B_2$ is uniquely defined by the following conditions. (1) $l_0$ coincides with $l^2$. (2) The endpoints of the arc $l_i$,\, $i = 0,\, 1,\, 2,\ldots$ are denoted by $\lam_i$ and $\nu_i$; $\lam_i$ lies on the segment $LB_2$ and $\nu_i$ lies on the segment $\nu_0B_2$. (3) $\lam_{i+1}$ lies on the segment $\nu_iA_2$.

Similarly, the sequence $c_0,\, c_1,\, c_2, \ldots$ of arcs of hyperbolas with foci at $A_2$ and $D_1$ is uniquely defined by the following conditions. (1) $c_0$ coincides with $c_K^1$. (2) The endpoints of the arc $c_i$,\, $i = 0,\, 1,\, 2,\ldots$ are denoted by $\s_i$ and $\chi_i$; $\s_i$ lies on the segment $C_1D_1$ and $\chi_i$ lies on the segment $\chi_0D_1$. (3) $\s_{i+1}$ lies on the extension of the segment $\chi_iA_2$.

Below we will need the following statement.

\begin{lem}\label{l1}
(a) Consider two different points $F_1$ and $F_2$, and let two rays from $F_1$ intersect two rays from $F_2$ at four points: $e_1$,\, $e_2$,\, $h_1$,\, $h_2$. Draw two ellipses $\mathcal{E}_1$,\, $\mathcal{E}_2$ and two hyperbolas $\mathcal{H}_1$,\, $\mathcal{H}_2$ with foci $F_1$ and $F_2$ through $e_1$,\, $e_2$,\, $h_1$,\, $h_2$, respectively. (Recall that by hyperbola we mean the corresponding branch of hyperbola.) We claim that the point of intersection of  $\mathcal{E}_1$ and $\mathcal{E}_2$, the point of intersection of  $\mathcal{H}_1$ and $\mathcal{H}_2$, and $F_2$ are collinear.

(b) Consider a ray from $F_2$ and two points $u_1$,\, $u_2$ on it. Draw the ellipse $\mathcal{E}_1$ and the hyperbola $\mathcal{H}_1$ with foci $F_1$,\, $F_2$ through $u_1$, and the ellipse $\mathcal{E}_2$ and the hyperbola $\mathcal{H}_2$ with the same foci $F_1$,\, $F_2$ through $u_2$. Take a ray from $F_1$ intersecting $\mathcal{E}_1$ and $\mathcal{H}_1$ at $e_1$ and $h_1$, respectively. Let the rays $F_2e_1$ and $F_2h_1$ intersect $\mathcal{E}_2$ and $\mathcal{H}_2$ at $e_2$ and $h_2$. We claim that the points $e_2$,\, $h_2$, and $F_1$ are collinear.
\end{lem}


\begin{figure}
\begin{picture}(0,150)
\rput(3.5,0){
\psline[linewidth=0.4pt](0,0)(0,4)
\psline[linewidth=0.4pt](0,0)(-2,5)
\psline[linewidth=0.4pt](2,3)(-2,5)
\psline[linewidth=0.4pt](2,3)(-0.91,2.27)
\psdots(0,0)(2,3)(-0.91,2.27)(0,2.5)(-2,5)(0,4)
\psline[linewidth=0.2pt,linestyle=dashed,linecolor=blue](2,3)(-2.5,3.45)
\pscurve[linewidth=0.8pt,linecolor=brown](-1.73,5.5)(-2,5)(-2.32,4.225)(-2.5,3.45)(-2.55,3.1)
\pscurve[linewidth=0.8pt,linecolor=brown](-0.78,2.2)(-0.91,2.27)(-1.8,2.85)(-2.5,3.45)(-2.8,3.8)
\pscurve[linewidth=0.8pt,linecolor=brown](0.1,2.4)(0,2.5)(-0.29,2.875)(-0.5,3.25)(-0.52,3.33)
\pscurve[linewidth=0.8pt,linecolor=brown](-0.53,3.15)(-0.5,3.25)(-0.29,3.625)(0,4)(0.23,4.2)
\rput(0.4,0){$F_1$}
\rput(2.4,3){$F_2$}
\rput(-2.6,4.5){$\mathcal{E}_1$}
\rput(-0.5,3.8){\scalebox{0.8}{$\mathcal{E}_2$}}
\rput(-2,2.6){$\mathcal{H}_1$}
\rput(-0.5,2.75){\scalebox{0.8}{$\mathcal{H}_2$}}
\rput(-1.7,5.1){\scalebox{0.9}{$e_1$}}
\rput(-0.1,4.3){\scalebox{0.9}{$e_2$}}
\rput(-1,2.05){\scalebox{0.9}{$h_1$}}
\rput(0.3,2.25){\scalebox{0.9}{$h_2$}}
\rput(-2.5,0.3){{\scalebox{1.2}{(a)}}}
\rput(9.5,0){
\psline[linewidth=0.2pt,linestyle=dashed,linecolor=blue](0,0)(0,4)
\psline[linewidth=0.4pt](0,0)(-2,5)
\psline[linewidth=0.4pt](2,3)(-2,5)
\psline[linewidth=0.4pt](2,3)(-0.91,2.27)
\psdots(0,0)(2,3)(-0.91,2.27)(-2,5)(-2.5,3.45)(-0.5,3.25)
\rput(-2.8,3.4){\scalebox{0.9}{$u_1$}}
\rput(-0.67,3.5){\scalebox{0.9}{$u_2$}}
\psline[linewidth=0.4pt](2,3)(-2.5,3.45)
\pscurve[linewidth=0.8pt,linecolor=brown](-1.73,5.5)(-2,5)(-2.32,4.225)(-2.5,3.45)(-2.55,3.1)
\pscurve[linewidth=0.8pt,linecolor=brown](-0.78,2.2)(-0.91,2.27)(-1.8,2.85)(-2.5,3.45)(-2.8,3.8)
\pscurve[linewidth=0.8pt,linecolor=brown](0.1,2.4)(0,2.5)(-0.29,2.875)(-0.5,3.25)(-0.52,3.33)
\pscurve[linewidth=0.8pt,linecolor=brown](-0.53,3.15)(-0.5,3.25)(-0.29,3.625)(0,4)(0.23,4.2)
\rput(0.4,0){$F_1$}
\rput(2.4,3){$F_2$}
\rput(-2.6,4.5){$\mathcal{E}_1$}
\rput(-0.01,3.6){\scalebox{0.8}{$\mathcal{E}_2$}}
\rput(-2,2.6){$\mathcal{H}_1$}
\rput(-0.5,2.75){\scalebox{0.8}{$\mathcal{H}_2$}}
\rput(-1.7,5.1){\scalebox{0.9}{$e_1$}}
\rput(-0.1,4.3){\scalebox{0.9}{$e_2$}}
\rput(-1,2.05){\scalebox{0.9}{$h_1$}}
\rput(0.3,2.25){\scalebox{0.9}{$h_2$}}
\rput(-2.5,0.3){{\scalebox{1.2}{(b)}}}
}
}
\end{picture}
\caption{}
\label{fig2}
\end{figure}


Its proof is given in Appendix. Note that another proof of Lemma using elegant geometric argument is provided by Pavel Dolgirev in \cite{dol}.

Let us now prove by induction that the points $\chi_i,\, \nu_i$, and $A_2$ (and therefore the points $\s_{i+1}$,\, $\chi_i$,\, $\nu_i$,\, $\lam_{i+1}$,\, $A_2$) are collinear. For $i=0$ this statement is obvious, since the points $\chi_0,\, \nu_0,\, A_2$ lie on the line $A_2H_1$. Now assume that $\chi_{i-1}$,\, $\nu_{i-1}$, and $A_2$ are collinear for $i \ge 1$, and prove that $\chi_{i}$,\, $\nu_{i}$,\, $A_2$ are collinear. Consider the dilation with the center at $A_2$ that takes $D_1$ to $B_2$. This dilation takes $c_i$ to the arcs $c_i'$ of hyperbolas with foci at  $A_2$ and $B_2$. Let the endpoints of $c_i'$ be $\s_i'$ and $\chi_i'$; the assumption of induction implies that the points $\s_{i}'$,\, $\chi_{i-1}'$,\, $\nu_{i-1}$,\, $\lam_{i}$,\, $A_2$ are collinear. We are going to prove that $\chi_i',\, \nu_i$, and $A_2$ are collinear.

The rays $A_2\s_{i-1}$ and $A_2\s_i$ intersect the rays $B_2A_1$ and $B_2\s_{i-1}'$ at four points: $\lam_{i-1}$,\, $\lam_i$,\, $\s_{i-1}'$,\, $\s_i'$. The ellipses $l_{i-1}$ and $l_i$ contain the points $\lam_{i-1}$ and $\lam_i$, and the hyperbolas $c_{i-1}'$ and $c_i'$ contain the points $\s_{i-1}'$, and $\s_i'$, therefore, according to the statement (a) of Lemma \ref{l1}, the point of intersection of $l_{i-1}$ and $c_{i-1}'$, the point of intersection of $l_i$ and $c_i'$, and $B_2$ are collinear. Further, using additionally that the points $\nu_{i-1}$ and $\nu_i$ lie on a single ray from $B_2$, $\chi_{i-1}'$ and $\chi_i'$ lie on another ray from $B_2$, $\nu_{i-1}$ and $\chi_{i-1}'$ lie on a ray from $A_2$, and applying the statement (b) of Lemma, we finally conclude that the points  $\chi_i',\, \nu_i$, and $A_2$ are collinear.

\begin{figure}
\begin{picture}(0,265)

\rput(10,8){
\scalebox{0.7}{
\psdots(-12,0)(0,-12)
\psdots[dotsize=2.5pt](0,-4)(0,4)
\psline[linewidth=0.4pt](-12,0)(6,6)
\psline[linewidth=0.4pt](-12,0)(2.4,-4.8)
\psline[linewidth=0.4pt](0,-12)(6,6)
\psline[linewidth=0.4pt](0,-12)(-4.8,2.4)

\psdots[dotsize=0.5pt](-1.5,-1.98)
(-2.7,-2.723)(-2.4,-2.485)(-2.1,-2.29)(-1.8,-2.12)(-1.5,-1.98)
(-1.26,-3.58)
(-0.561,-3.813)
(-0.603,-3.188)
(-0.26,-3.65)
(-1.04,-3.413)(-0.82,-3.282)
(-0.46,-3.741)(-0.36,-3.688)
\psline[linewidth=0.1pt,linestyle=dashed,linecolor=brown](0,-4)(-1.5,-1.98)
\psline[linewidth=0.1pt,linestyle=dashed,linecolor=brown](0,-12)(-2.4,4.032) 
\psline[linewidth=0.1pt,linestyle=dashed,linecolor=brown](0,-12)(-1.1457,4.7428) 
\psline[linewidth=0.1pt,linestyle=dashed,linecolor=brown](0,-12)(-0.52,4.7) 
\pscurve[linecolor=blue,linewidth=0.4pt](-3,-3)(-2.7,-2.723)(-2.4,-2.485)(-2.1,-2.29)(-1.8,-2.12)(-1.5,-1.98) 
\pscurve[linecolor=blue,linewidth=0.4pt](-1.26,-3.58)(-1.04,-3.413)(-0.82,-3.282)(-0.603,-3.188) 
\pscurve[linecolor=blue,linewidth=0.4pt](-0.561,-3.813)(-0.46,-3.741)(-0.36,-3.688)(-0.26,-3.65) 
\psdots[dotsize=0.5pt](-1.99,1.325)
(-1.07,3.644)
\psecurve[linecolor=blue,linewidth=0.4pt](-1.2476,1.1611)(-1.99,1.32)(-3.229,1.7053)(-4.8,2.4)(-7.2824,3.8074) 
\psline[linewidth=0.1pt,linestyle=dashed,linecolor=brown](0,4)(-1.99,1.32)
\psecurve[linecolor=blue,linewidth=0.4pt](-0.4667,2.5356)(-1,2.653)(-1.41885,2.7979)(-2.28,3.24)(-2.4279,3.3331) 
\psecurve[linecolor=blue,linewidth=0.4pt](0,3.2764)(-0.48,3.367)(-0.6694,3.4224)(-1.07,3.644)(-1.3656,3.8663) 
\psdots[dotsize=0.7pt](-0.857,-3.3)(-1.46,2.815)(-2.943,1.605)


\rput(-12.1,-0.5){\scalebox{1.4}{$A_1$}}
\rput(0.5,-12.2){\scalebox{1.4}{$A_2$}}
\rput(0.2,-4.5){\scalebox{1.2}{$B_2$}}
\rput(-0.2,4.4){\scalebox{1.2}{$D_1$}}
\rput(-3.5,1.5){\scalebox{1.2}{$c_0$}}
\rput(-4.9,2.7){\scalebox{1.2}{$\s_0$}}
\rput(-2.5,3.4){\scalebox{1}{$\s_1$}}
\rput(-2.2,1.2){\scalebox{1}{$\chi_0$}}
\rput(-1.7,2.7){\scalebox{1}{$c_1$}}
\rput(-0.8,2.5){\scalebox{1}{$\chi_1$}}
\rput(-2.5,-2.25){\scalebox{1}{$l_0$}}
\rput(-1,-3.05){\scalebox{1}{$l_1$}}
\rput(-3.2,-3.2){\scalebox{1}{$\lam_0$}}
\rput(-1.25,-1.85){\scalebox{1}{$\nu_0$}}
\rput(-1.5,-3.7){\scalebox{0.9}{$\lam_1$}}
\rput(-0.35,-3.1){\scalebox{1}{$\nu_1$}}

\rput(-2.7,-4.8){\scalebox{1}{$\s_0'$}}
\rput(-0.95,-5.5){\scalebox{1}{$\chi_0'$}}
\rput(-1.8,-5.4){\scalebox{1}{$c_0'$}}
  \rput(-1.2,-4.1){\scalebox{0.9}{$\s_1'$}}
  \rput(-0.35,-4.85){\scalebox{0.9}{$\chi_1'$}}
  \rput(-0.9,-4.75){\scalebox{0.9}{$c_1'$}}

\rput(-0.1,-6){
\scalebox{0.5}{
\psecurve[linecolor=blue,linewidth=0.8pt](-1.2476,1.1611)(-1.99,1.32)(-3.229,1.7053)(-4.8,2.4)(-7.2824,3.8074) 
\psline[linewidth=0.1pt,linestyle=dashed,linecolor=brown](0,4)(-1.99,1.32)
\psecurve[linecolor=blue,linewidth=0.8pt](-0.4667,2.5356)(-1,2.653)(-1.41885,2.7979)(-2.28,3.24)(-2.4279,3.3331) 
\psline[linewidth=0.4pt](-4.8,2.4)(0,4)
\psline[linewidth=0.1pt,linestyle=dashed,linecolor=brown](0,4)(-1.99,1.32)
}}
}}
\end{picture}
\caption{}
\label{fig3}
\end{figure}

The resulting body is the union of sets $Q_L$,\, $Q_{C_1}$,\, $Q_{C_2}$,\, $Q_K$, and 8 sequences of arcs. Let us show that it is invisible from the points $A_1$ and $A_2$.

It is enough to show invisibility from the point $A_2$; the invisibility from $A_1$ can be verified in a completely similar way. If the light ray emanating from $A_2$ does not belong to the angle $C_1A_2K$ (which is the union of the angles $C_1A_2O$ and $OA_2K$), it does not hit the body. It remains to consider the cases when it belongs to the angle $C_1A_2O$ and to the angle $OA_2K$. We consider only the angle $C_1A_2O$, since the case of the angle $OA_2K$ is completely similar.

The light ray reflects from an arc $l_{i-1}$ ($i \ge 1$) and then goes along a straight line containing $B_2$. (In Fig. \ref{fig4} the case $i=1$ is shown.) Next it reflects from $l_i$ and then goes along a line containing $A_2$. Then after reflection from the arc $c_i$ it goes along a line containing $D_1$, hits $c_{i-1}$ and finally goes along a line containing $A_2$. (The points of reflection are denoted by $P_1$,\, $P_2$,\, $P_3$,\, $P_4$.) It remains to prove that this final line coincides with the initial one.

\begin{figure}
\begin{picture}(0,265)

\rput(10,8){
\scalebox{0.7}{
\psdots(-12,0)(0,-12)
\psdots[dotsize=2.5pt](0,-4)(0,4)
\psline[linewidth=0.4pt](-12,0)(6,6)
\psline[linewidth=0.4pt](-12,0)(2.4,-4.8)
\psline[linewidth=0.4pt](0,-12)(6,6)
\psline[linewidth=0.4pt](0,-12)(-4.8,2.4)

\psdots[dotsize=0.5pt](-1.5,-1.98)
(-2.7,-2.723)(-2.4,-2.485)(-2.1,-2.29)(-1.8,-2.12)(-1.5,-1.98)
(-1.26,-3.58)
(-0.561,-3.813)
(-0.603,-3.188)
(-0.26,-3.65)
(-1.04,-3.413)(-0.82,-3.282)
(-0.46,-3.741)(-0.36,-3.688)
\psline[linewidth=0.1pt,linestyle=dashed,linecolor=brown](0,-4)(-1.5,-1.98)
\psline[linewidth=0.1pt,linestyle=dashed,linecolor=brown](0,-12)(-2.4,4.032) 
\psline[linewidth=0.1pt,linestyle=dashed,linecolor=brown](0,-12)(-1.1457,4.7428) 
\psline[linewidth=0.1pt,linestyle=dashed,linecolor=brown](0,-12)(-0.52,4.7) 
\pscurve[linecolor=blue,linewidth=0.4pt](-3,-3)(-2.7,-2.723)(-2.4,-2.485)(-2.1,-2.29)(-1.8,-2.12)(-1.5,-1.98) 
\pscurve[linecolor=blue,linewidth=0.4pt](-1.26,-3.58)(-1.04,-3.413)(-0.82,-3.282)(-0.603,-3.188) 
\pscurve[linecolor=blue,linewidth=0.4pt](-0.561,-3.813)(-0.46,-3.741)(-0.36,-3.688)(-0.26,-3.65) 
\psdots[dotsize=0.5pt](-1.99,1.325)
(-1.07,3.644)
\psecurve[linecolor=blue,linewidth=0.4pt](-1.2476,1.1611)(-1.99,1.32)(-3.229,1.7053)(-4.8,2.4)(-7.2824,3.8074) 
\psline[linewidth=0.1pt,linestyle=dashed,linecolor=brown](0,4)(-1.99,1.32)
\psecurve[linecolor=blue,linewidth=0.4pt](-0.4667,2.5356)(-1,2.653)(-1.41885,2.7979)(-2.28,3.24)(-2.4279,3.3331) 
\psecurve[linecolor=blue,linewidth=0.4pt](0,3.2764)(-0.48,3.367)(-0.6694,3.4224)(-1.07,3.644)(-1.3656,3.8663) 
\psdots[dotsize=0.7pt](-0.857,-3.3)(-1.46,2.815)(-2.943,1.605)

\psline[linecolor=red,linewidth=0.4pt,arrows=->,arrowscale=1.7](0,-12)(-2.1,-2.29)
\psline[linecolor=red,linewidth=0.4pt,arrows=->,arrowscale=1.7](-2.1,-2.29)(-0.857,-3.3)(-1.46,2.815)
\psline[linecolor=red,linewidth=0.4pt,arrows=->,arrowscale=1.7](-1.46,2.815)(-2.943,1.605)(-4.2,7.42)

\psline[linecolor=red,linestyle=dashed,linewidth=0.4pt](-0.857,-3.3)(-0.73,-4.585)(-1.47,-5.195)

\rput(-12.1,-0.5){\scalebox{1.4}{$A_1$}}
\rput(0.5,-12.2){\scalebox{1.4}{$A_2$}}

\rput(0.1,-4.45){\scalebox{1.2}{$B_2$}}
\rput(-0.2,4.4){\scalebox{1.2}{$D_1$}}

\rput(-4.2,1.85){\scalebox{0.8}{$c_{i-1}$}}
\rput(-3,1.3){\scalebox{1}{$P_{4}$}}
\rput(-2,2.85){\scalebox{0.9}{$c_i$}}
\rput(-1.4,3.1){\scalebox{1}{$P_{3}$}}
\rput(-2.85,-2.45){\scalebox{0.8}{$l_{i-1}$}}
\rput(-2.1,-1.95){\scalebox{1}{$P_{1}$}}
\rput(-1.35,-3.3){\scalebox{0.9}{$l_i$}}
\rput(-0.7,-3.5){\scalebox{1}{$P_{2}$}}
\rput(-2,-5.12){\scalebox{0.8}{$c_{i-1}'$}}
\rput(-1.6,-5.5){\scalebox{0.8}{$P_{4}'$}}
  \rput(-1.2,-4.65){\scalebox{0.9}{$c_i'$}}
  \rput(-0.7,-4.9){\scalebox{0.8}{$P_{3}'$}}

\rput(-0.1,-6){
\scalebox{0.5}{
\psecurve[linecolor=blue,linewidth=0.8pt](-1.2476,1.1611)(-1.99,1.32)(-3.229,1.7053)(-4.8,2.4)(-7.2824,3.8074) 
\psline[linewidth=0.1pt,linestyle=dashed,linecolor=brown](0,4)(-1.99,1.32)
\psecurve[linecolor=blue,linewidth=0.8pt](-0.4667,2.5356)(-1,2.653)(-1.41885,2.7979)(-2.28,3.24)(-2.4279,3.3331) 
\psline[linewidth=0.4pt](-4.8,2.4)(0,4)
\psline[linewidth=0.1pt,linestyle=dashed,linecolor=brown](0,4)(-1.99,1.32)
}}
}}
\end{picture}
\caption{}
\label{fig4}
\end{figure}

Consider the arcs $c_{i-1}'$,\, $c_i'$ homothetical to $c_{i-1}$,\, $c_i$ and extend the segment $P_2P_3$ of trajectory between $l_i$ and $c_i$ until the intersection $P_3'$ with $c_i'$ (dashed red line in the figure). The point $P_3'$ is homothetic to $P_3$. Next we take the straight line through $P_3'$ and $B_2$ and fix the point $P_4'$ of its intersection with $c_{i-1}'$. The segment $P_3'P_4'$ is a dashed red line in the figure. The point $P_4'$ is homothetic to $P_4$. We already know that the point of intersection of the hyperbola $c_{i-1}'$ with the ellipse $l_{i-1}$, the point of intersection of the hyperbola $c_{i}'$ with the ellipse $l_{i}$, and $B_2$ are collinear. The points $P_2$,\, $P_3'$, and $A_2$ are collinear, the points $P_4'$,\, $P_3'$, and $B_2$ are collinear, and the points $P_1$,\, $P_2$, and $B_2$ are also collinear; therefore, according to statement (b) of Lemma \ref{l1}, the points $P_1'$,\, $P_4'$, and $A_2$ are also collinear. Thus, $P_4'$ lies on $A_2P_1$, and therefore, $P_4$ also lies on this line. This implies that the line $A_2P_4$ of final motion of the particle coincides with the line $A_2P_1$ of initial motion. Invisibility is proved.

\begin{figure}[h]
\centering
\includegraphics[scale = 1]{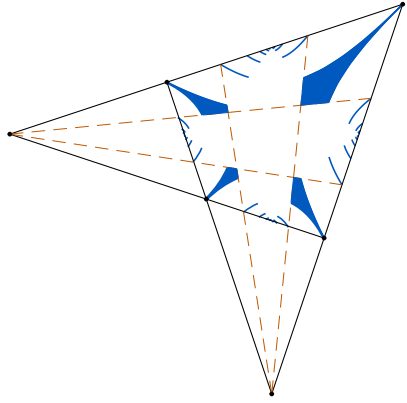}
\caption{Body invisible from two points, 2D case}
\label{fig:v03}
\end{figure}

The invisible body is shown in Fig. \ref{fig:v03}.
It is the union of four curvilinear quadrangles and 8 infinite sequences of curves of vanishing length.

\begin{figure}[h]
\centering
\includegraphics[width=400pt]{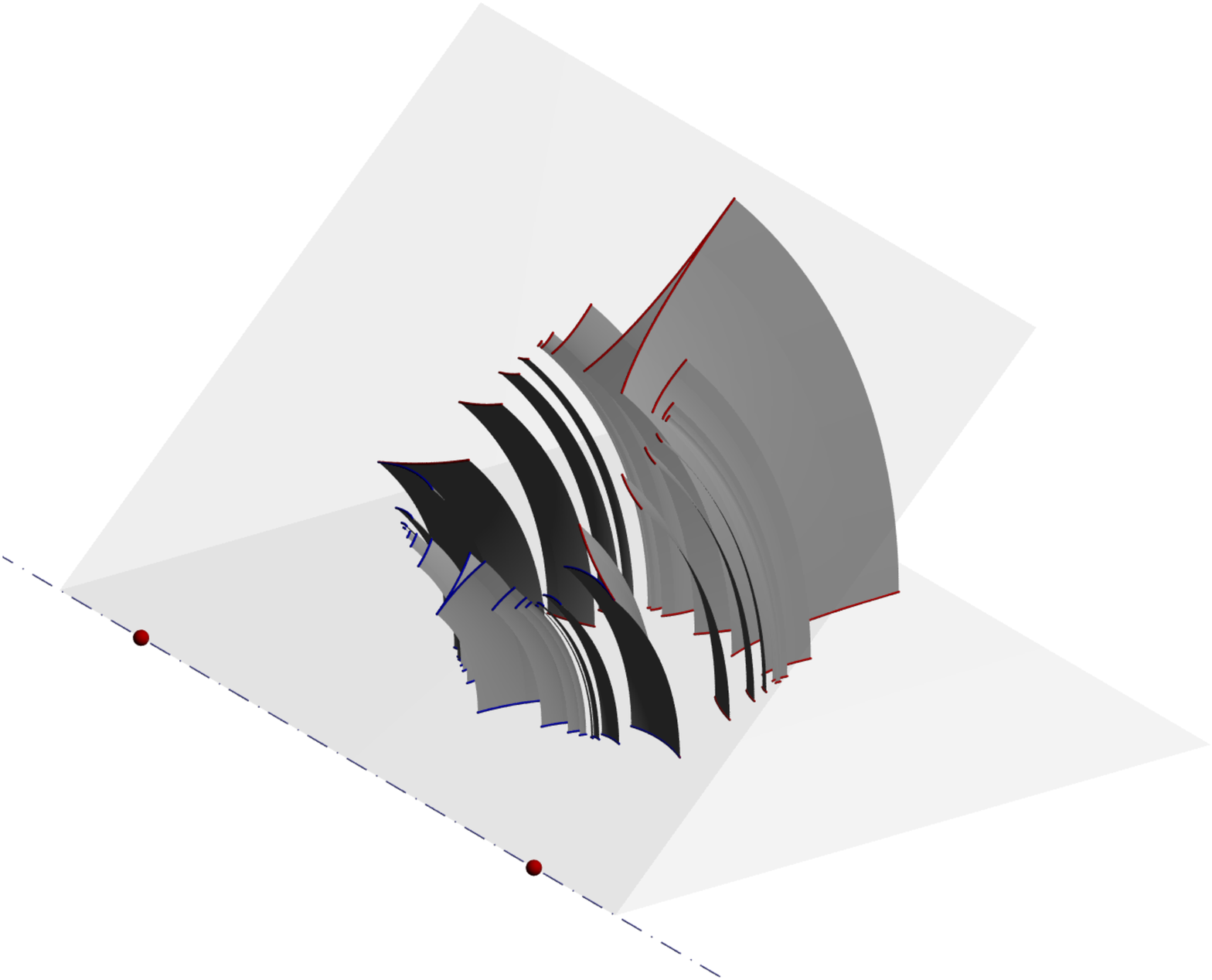}\\
\caption{Body invisible from two points, 3D case}
\label{fig:v04}
\end{figure}

To obtain a body invisible from two points in a higher dimensional settings, it is sufficient to rotate the two-dimensional construction around the axis $A_1A_2$ (see Fig.~\ref{fig:v04}
for the illustration of the 3D case). It is not necessary to make a full rotation; for instance, in the 3D case a rotation by a certain angle (not necessary $360^0$) results in an invisible body. We hence proved Theorem~\ref{t1}.

\begin{zam}
The construction of the invisible body does not have to be symmetric. In fact, it is enough to choose four points $H_1,\, M,\, H_2,\, N$ in the quadrangles $OB_1C_1D_1$,\, $OD_1KD_2$,\, $OD_2C_2B_2$,\, $OB_2LB_1$, respectively. We require that the straight lines $H_1M$ and $H_2N$ contain $A_1$, the lines $MH_2$ and $NH_1$ contain $A_2$, and each of the points $H_1,\, M,\, H_2,\, N$ belongs to the exterior of both curves through $C_1,\, K,\, C_2$, or $L$, respectively. Moreover, we require that each of the segments $A_1N$,\, $A_2N$ intersects only one curve through $L$ and that the analogous conditions related to the other three vertices of $LC_1KC_2$ hold.
\end{zam}

To demonstrate the existence of a body invisible in one direction and from one point simultaneously, a similar construction may be used. The difference is that here one needs to use 4 parabolas, 2 ellipses, and 2 hyperbolas instead of 4 ellipses and 4 hyperbolas. The details of this modified construction are left to the reader.

\section*{Appendix}

Consider two rays emanating from $F_2$ and a ray emanating from $F_1$. Let the points of intersection of the first two ones with the third one be denoted by $h$ and $e$. Let $\al = \measuredangle F_1F_2h$,\, $\bt = \measuredangle F_1F_2e$,\, $\gam = \measuredangle hF_1F_2$. Draw an ellipse and a hyperbola with foci at $F_1$ and $F_2$ through $e$ and $h$, respectively; let $u$ be the point of their intersection and $\vphi = \measuredangle F_1F_2u$. We are going to show that $\vphi$ depends only on $\al$ and $\bt$, and does not depend on $\gam$ (that is, does not depend on the choice of the ray from $F_1$). Thus the statement (a) of Lemma \ref{l1} will be proved. The proof of the statement (b) is completely similar and therefore is omitted here.

\begin{figure}
\begin{picture}(0,130)

\rput(8,0){
\psline[linewidth=0.2pt,linestyle=dashed,linecolor=blue](0,0)(2,3)
\psline[linewidth=0.4pt](0,0)(-2,5)
\psline[linewidth=0.4pt](2,3)(-2,5)
\psline[linewidth=0.4pt](2,3)(-0.91,2.27)
\psdots(0,0)(2,3)(-0.91,2.27)(-2,5)(-2.5,3.45)
\rput(-2.8,3.4){\scalebox{0.9}{$u$}}
\psline[linewidth=0.4pt](2,3)(-2.5,3.45)
\pscurve[linewidth=0.8pt,linecolor=brown](-1.73,5.5)(-2,5)(-2.32,4.225)(-2.5,3.45)(-2.55,3.1)
\pscurve[linewidth=0.8pt,linecolor=brown](-0.78,2.2)(-0.91,2.27)(-1.8,2.85)(-2.5,3.45)(-2.8,3.8)
\rput(0.4,0){$F_1$}
\rput(2.4,3){$F_2$}
\rput(-2.6,4.5){$\mathcal{E}$}
\rput(-2,2.6){$\mathcal{H}$}
\rput(-1.7,5.1){\scalebox{0.9}{$e$}}
\rput(-1,2.05){\scalebox{0.9}{$h$}}
\psarc[linecolor=blue](2,3){0.5}{174}{236}
\psarc[linecolor=blue](2,3){0.45}{174}{236}
\rput(1.85,2.45){$\vphi$}
\psarc[linecolor=blue](2,3){0.25}{153.3}{236}
\psarc[linecolor=blue](2,3){0.275}{153.3}{236}
\psarc[linecolor=blue](2,3){0.3}{153.3}{236}
\rput(1.8,3.35){$\bt$}
\psarc[linecolor=blue](2,3){0.7}{194}{236}
\rput(1.33,2.4){$\al$}
\psarc[linecolor=blue](0,0){0.3}{56.3}{112}
\rput(0,0.5){$\gam$}
\rput(0.2,2.3){$b_1$}
\rput(0,4.3){$b_2$}
\rput(-0.7,1.1){$a_1$}
\rput(-1.3,4){$a_2$}
\rput(1.15,1.35){$f$}
\rput(-0.2,3.4){$c$}
}
\end{picture}
\caption{}
\label{fig5}
\end{figure}

Denote $a_1 = F_1h$,\, $a_2 = F_1e$,\, $b_1 = F_2h$,\, $b_2 = F_2e$,\, $f = F_1F_2$,\, $c = F_2u$. By the focal property of hyperbola and ellipse we have
$$
F_1h - F_2h = F_1u - F_2u \quad \text{and} \quad F_1e + F_2e = F_1u + F_2u,
$$
whence
$$
a_1 - b_1 = \sqrt{f^2 + c^2 - 2cf\cos\vphi} - c, \quad\quad  a_1 + b_1 = \sqrt{f^2 + c^2 - 2cf\cos\vphi} + c,
$$
and therefore,
\beq\label{1}
\sqrt{f^2 + c^2 - 2cf\cos\vphi} = \frac 12 (a_1 - b_1 + a_2 + b_2),
\eeq
\beq\label{2}
c = \frac 12 (-a_1 + b_1 + a_2 + b_2).
\eeq
Further, by sine law,
$$
\frac{a_1}{\sin\al} = \frac{b_1}{\sin\gam} = \frac{f}{\sin(\al+\gam)}, \quad\quad \frac{a_2}{\sin\bt} = \frac{b_2}{\sin\gam} = \frac{f}{\sin(\bt+\gam)},
$$
hence
$$
a_1 = f\, \frac{\sin\al}{\sin(\al+\gam)}, \quad b_1 = f\, \frac{\sin\gam}{\sin(\al+\gam)}, \quad a_2 = f\, \frac{\sin\bt}{\sin(\bt+\gam)}, \quad b_2 = f\, \frac{\sin\gam}{\sin(\bt+\gam)},
$$
and substituting these quantities in (\ref{1}) and (\ref{2}) one gets
$$
f^2 + c^2 - 2cf\cos\vphi = \frac{f^2}{4}\, \Big[ \frac{\sin\al-\sin\gam}{\sin(\al+\gam)} + \frac{\sin\bt+\sin\gam}{\sin(\bt+\gam)} \Big]^2,
$$
$$
c = \frac{f}{2}\, \Big[ \frac{\sin\gam-\sin\al}{\sin(\al+\gam)} + \frac{\sin\bt+\sin\gam}{\sin(\bt+\gam)} \Big],
$$
and thus,
$$
\cos\vphi = \frac AB,
$$
where
$$
B = 2cf = f^2\, \Big[ \frac{\sin\gam-\sin\al}{\sin(\al+\gam)} + \frac{\sin\bt+\sin\gam}{\sin(\bt+\gam)} \Big],
$$
$$
A = f^2 + c^2 - \frac{f^2}{4}\, \Big[ \frac{\sin\al-\sin\gam}{\sin(\al+\gam)} + \frac{\sin\bt+\sin\gam}{\sin(\bt+\gam)} \Big]^2 =
$$
$$
= f^2 \left\{ 1 + \frac 14\, \Big[ \frac{\sin\gam-\sin\al}{\sin(\al+\gam)} + \frac{\sin\bt+\sin\gam}{\sin(\bt+\gam)} \Big]^2 - \frac{1}{4}\, \Big[ \frac{\sin\al-\sin\gam}{\sin(\al+\gam)} + \frac{\sin\bt+\sin\gam}{\sin(\bt+\gam)} \Big]^2 \right\}.
$$
Thus,
$$
B\, \frac{\sin(\al+\gam)\sin(\bt+\gam)}{f^2} = \sin(\bt+\gam)(\sin\gam-\sin\al) + \sin(\al+\gam)(\sin\bt+\sin\gam),
$$
$$
A\, \frac{\sin(\al+\gam)\sin(\bt+\gam)}{f^2} = \sin(\al+\gam)\sin(\bt+\gam) + (\sin\gam-\sin\al)(\sin\gam+\sin\bt).
$$

After  some algebra one gets
$$
B\, \frac{\sin(\al+\gam)\sin(\bt+\gam)}{f^2} = 2\sin\gam \cos\frac{\al-\bt}{2}\, \Big[ \sin\Big( \gam + \frac{\al+\bt}{2} \Big) + \sin\frac{\bt-\al}{2} \Big],
$$
$$
A\, \frac{\sin(\al+\gam)\sin(\bt+\gam)}{f^2} = 2\sin\gam \cos\frac{\al+\bt}{2}\, \Big[ \sin\Big( \gam + \frac{\al+\bt}{2} \Big) + \sin\frac{\bt-\al}{2} \Big],
$$
and therefore,
$$
\cos\vphi = \frac{\cos\frac{\al+\bt}{2}}{\cos\frac{\al-\bt}{2}}.
$$

\section*{Acknowledgements}

This work was supported in part by FEDER funds through COMPETE – Operational Programme Factors of Competitiveness ("Programa Operacional Factores de Competitividade") and by Portuguese funds through the Center for Research and Development in Mathematics and Applications and the Portuguese Foundation for Science and Technology ("FCT–Fundacao para a Ciencia e a Tecnologia"), within project PEst-C/MAT/UI4106/2011 with COMPETE number FCOMP-01-0124-FEDER-022690, as well as by the FCT research project PTDC/MAT/113470/2009.

\end{document}